\documentclass[12pt]{amsart} 
\usepackage{amsthm, amscd}
\usepackage{amssymb,latexsym} 
\input xy
\xyoption{all}

\theoremstyle{plain}
\newtheorem*{introtheorem}{Theorem}
\newtheorem*{introcorollary}{Corollary}
\newtheorem{theorem}{Theorem}[section]
\newtheorem{proposition}[theorem]{Proposition}
\newtheorem{corollary}[theorem]{Corollary}
\newtheorem{lemma}[theorem]{Lemma}
 
\theoremstyle{definition}    
\newtheorem{definition}[theorem]{Definition}
\newtheorem{remark}[theorem]{Remark}
\newtheorem{example}[theorem]{Example}
   
\newtheorem{conjecture}[theorem]{Conjecture}
\newtheorem{conjecture/question}[theorem]{Conjecture/Question} 
\newtheorem{question}[theorem]{Question}

\theoremstyle{remark}

\newcommand{\PP}{\mathbf{P}}
\newcommand{\ZZ}{\mathbf{Z}}
\newcommand{\CC}{\mathbf{C}}

\newcommand{\OO}{\mathcal{O}}

\newcommand{\LL}{\mathcal{L}}

\newcommand{\dt}{\text{det}}

\begin{document}

\title{\bf Dimension estimates for Hilbert schemes and effective base point freeness 
on moduli spaces of vector bundles on curves}
\author{Mihnea Popa}
\address{Department of Mathematics
\\University of Michigan \\ 525 East University \\ Ann Arbor, MI 48109-1109}
\address{Institute of
Mathematics of the Romanian Academy}
\email{mpopa@math.lsa.umich.edu}

\maketitle
\markboth{MIHNEA POPA}{\bf HILBERT SCHEMES AND EFFECTIVE BASE POINT FREENESS}

$$\textbf{Introduction}$$

It is a well established fact that the solutions of many problems involving families of 
vector bundles should essentially depend on good estimates for the dimension of the Hilbert schemes of 
coherent quotients of a given bundle. Deformation theory provides basic cohomological 
dimension bounds, but most of the time the cohomology groups involved are hard to 
estimate accurately and moreover do not provide bounds that work uniformly. On smooth algebraic 
curves, an optimal answer to this problem has been previously given only in the case of 
quotients of minimal degree by Mukai and Sakai. If $E$ is a vector bundle of rank $r$ and 
$$f_{k}=f_{k}(E):=\underset{{\rm rk}Q=k}{{\rm min}}\{{\rm deg}Q|~E\rightarrow Q\rightarrow 0\},$$
is the minimal degree of a quotient of $E$ of rank $k$, they show in \cite{Mukai} that
${\rm dim~}{\rm Quot}_{k,f_{k}}(E)\leq k(r-k)$, where in general  
${\rm Quot}_{k,d}(E)$ stands for Grothendieck's Hilbert (or Quot) scheme of coherent 
quotients of $E$ of rank $k$ and degree $d$.

In the first part of this paper we give an upper bound for the dimension of 
${\rm Quot}_{k,d}(E)$ for any degree $d$. The bound involves in an essential (and somewhat 
unexpected) way the invariant $f_{k}$. Examples provided in \ref{40} show that it is 
optimal at least in the case corresponding to line subbundles.

\begin{introtheorem}
Let $E$ be an arbitrary vector bundle of rank $r$ on a smooth projective curve $X$ over an 
algebraically closed field. Then:
$${\rm dim~}{\rm Quot}_{k,d}(E)\leq k(r-k)+(d-f_{k})(k+1)(r-k).$$
\end{introtheorem}
This generalizes (and uses) the result from \cite{Mukai} mentioned above, which is exactly 
the case $d=f_{k}$. The proof is based on 
induction on the difference $d-f_{k}$ and the key ingredient is a technique that allows one 
to ``eliminate'' all the minimal quotient bundles of $E$ while preserving a fixed nonminimal one.
This is achieved via elementary transformations along zero-dimensional subschemes of arbitrary length. 
The problem has also been previously given some precise answers in the particular case of generic
stable bundles in \cite{Teixidor} and \cite{Brambila}. In this case the dimension (and $f_{k}$) 
can be computed exactly (cf. Example \ref{40} below).

In the second part of the paper, we apply the estimate above to the basic problem of giving 
effective bounds for the global generation of multiples of generalized theta line bundles
on moduli spaces of vector bundles on curves.
On a smooth projective complex curve $X$ of genus $g\geq 2$, denote by $U_{X}(r,e)$
the moduli space of equivalence classes of semistable vector bundles on $X$ of rank $r$
and degree $e$ and by $SU_{X}(r,L)$ the moduli space of semistable rank $r$ vector bundles
of fixed determinant $L\in {\rm Pic}^{e}(X)$. When the choice of $L$ is of no special importance 
we will use the notation $SU_{X}(r,e)$ and furthermore
$SU_{X}(r)$ will denote the moduli space of vector bundles with trivial determinant.
It is well known (see \cite{Drezet}, Theorem $B$) that   
${\rm Pic}(SU_{X}(r,e))\cong \ZZ\cdot \LL$,
where the ample generator $\LL$ is called the determinant line bundle. Since the linear system $|\LL|$ is 
known to have base points in general, the main problem is to give effective 
bounds for the base point freeness of linear series of the form $|\LL^{p}|$
for some integer $p$. 
The particular shape of the Picard group implies that this is in fact a Fujita type problem for 
adjoint linear series. Previous work in this direction appears in papers of Le Potier \cite{Le Potier}
and Hein \cite{Hein}, where the authors obtain bounds with order of magnitude in the range 
suggested by Fujita's conjecture ($p>\frac{r^{2}}{4}(g-1)$ and $p> (r-1)^{2}(g-1)$ respectively). If 
$h$ denotes the greatest common divisor of $r$ and $e$, 
our main result is then formulated as follows: 

\begin{introtheorem}
$|\LL^{p}|$ is base point free on $SU_{X}(r,e)$ for $p\geq {\rm  max}\{
\frac{(r+1)^{2}}{4r} h,\frac{r^{2}}{4s}h\}$.
\end{introtheorem}
The number $s$ is an invariant associated to the space
$SU_{X}(r,e)$, which will be defined precisely in Section 4. We
will restrict here to saying that always $s\geq h$, so the bound
in the theorem is at most quadratic in the rank $r$, no worse
than $\frac{(r+1)^{2}}{4}$. The theorem substantially improves the existing bounds mentioned above. 
Most importantly, beyond the concrete numerology, the main improvement is that the present 
bound is independent of the genus of the curve, 
as it is natural to expect. The idea is again to make effective the method of Faltings
\cite{Faltings}, but the technique involved is substantially different from those used in 
\cite{Le Potier} or \cite{Hein}, the key ingredient being precisely the result 
on Hilbert schemes described above.

Two cases that have traditionally been under intensive study are 
that of degree $0$ bundles and that of bundles of degree
$\pm 1$
(or more generally $e\equiv 0({\rm mod}~r)$ and $e\equiv
\pm 1({\rm mod}~r))$. In the first  case
the bound that we obtain is quadratic in $r$, 
but somewhat surprisingly in the second case it is linear.

\begin{introcorollary}
(i) $|\LL^{p}|$ is base point free on $SU_{X}(r)$ for $p\geq \frac{(r+1)^{2}}{4}$.
\newline
(ii) $|\LL^{p}|$ is base point free on $SU_{X}(r,1)$ and
$SU_{X}(r,-1)$ for $p\geq r-1$.
\end{introcorollary}     
In fact in the first case one can do slightly better for $r$ even (see \S4).

Furthermore, making use of the notion of Verlinde bundle introduced in \cite{Popa2},
we find similar effective bounds for the base point freeness of
linear series on $U_{X}(r,e)$ (see also \cite{Popa2} \S5 and \S6).
The result (Theorem \ref{12}) can be formulated as follows:

\begin{introtheorem}
Let $F$ be a vector bundle of rank $\frac{r}{h}$ and degree $\frac{r}{h} (g-1)-\frac{e}{h}$ 
on $X$ and let $\Theta_{F}$ be the corresponding
generalized theta divisor on $U_{X}(r,e)$. Then $|p\Theta_{F}|$ is base point free for 
$$p\geq {\rm  max}\{
\frac{(r+1)^{2}}{4r} h,\frac{r^{2}}{4s}h\}.$$
\end{introtheorem}
In fact this statement is a special case of a result about linear series of a more 
general type (cf. Theorem \ref{11}).

In a different direction, we follow Le Potier's idea \cite{Le Potier} \S3 to observe that 
the bound given in the main theorem improves substantially the analogous bound for the 
global generation of multiples of the Donaldson determinant line bundle on moduli spaces  
of semistable sheaves on surfaces, the independence on the genus being again the crucial fact.
In particular this bounds the dimension of a projective space which is 
an ambient space for an embedding of the moduli space of $\mu$-semistable sheaves. 
In the case of rank $2$ (and degree $0$) sheaves this space 
in turn is known to be homeomorphic to the Donaldson-Uhlenbeck compactification of the 
moduli space of $ASD$-connections in gauge theory. 

\begin{introtheorem}
Let $(X,\OO_{X}(1))$ be a polarized smooth projective surface and $L$ a line bundle on $X$. 
Let $M=M_{X}(r,L,c_{2})$ be the 
moduli space of semistable sheaves of rank $r$, fixed determinant 
$L$ and second Chern class $c_{2}$ on $X$ and denote $n={\rm deg}(X)=\OO_{X}(1)^{2}$ and 
$d=n[\frac{r^{2}}{2}]$. If $\mathcal{D}$ is the Donaldson determinant line bundle
on $M$, then $\mathcal{D}^{\otimes p}$ is globally generated for $p\geq d\cdot \frac{(r+1)^{2}}{4}$
divisible by $d$.
\end{introtheorem}

The paper is organized as follows: in the first section we review a few basic facts 
about generalized theta divisors on moduli spaces of vector bundles in the context of our 
problem. 
In Section 2 we turn our attention to the dimension bounds for Hilbert schemes 
of coherent quotients of a given vector bundle. This paragraph is of a somewhat different flavor 
from the rest of the paper and can be read independently. 
Section 3 treats 
the special case of rank 2 vector bundles. We prove a well-known result of Raynaud
\cite{Raynaud} by a method intended to be a toy version of the proof of the general 
theorem in the subsequent section.
The fourth section contains the 
proof of the main base point freeness result on $SU_{X}(r,e)$, while the fifth treats the 
case of linear series on $U_{X}(r,e)$. There we also formulate some questions about optimal bounds 
in arbitrary rank, which for example in degree $0$ follow from our results in the case 
of rank $2$ and rank $3$ vector bundles. 
The last section is devoted to a brief treatment of the above mentioned application 
to moduli spaces of sheaves on surfaces.

\textbf{Acknowledgements.} I would especially like to thank my advisor, R. Lazarsfeld, whose 
support and suggestions have been decisive to this work. I am also indebted to I. Coand\u a, 
I. Dolgachev, W. Fulton, M. Roth and M. Teixidor I Bigas for very valuable discussions. 
In particular numerous conversations with M. Roth had a significant influence on \S2 below.

\section{\textbf{Background}}

The underlying idea for studying linear series on the moduli space $SU_{X}(r,e)$ has its
roots in the paper of Faltings \cite{Faltings}, where a construction of the moduli space 
based on theta divisors is given. A very nice introduction to the subject is provided 
in \cite{Beauville}. 

Fix $r$ and $e$ and denote $h={\rm gcd}(r,e)$, $r_{1}=\frac{r}{h}$ and $e_{1}=\frac{e}{h}$.
Consider a vector bundle $F$ of rank $pr_{1}$ and degree $p(r_{1}(g-1)-e_{1})$. Generically 
such a choice determines (cf. \cite{Drezet} 0.2) a \emph{theta divisor} $\Theta_{F}$ on $SU_{X}(r,e)$, 
supported on the set 
$$\Theta_{F}=\{E|~h^{0}(E\otimes F)\neq 0\}.$$

All the divisors $\Theta_{F}$ for $F\in U_{X}(pr_{1}, p(r_{1}(g-1)-e_{1}))$ belong to the 
linear system $|\LL^{p}|$, 
where $\LL$ is the \emph{determinant} line bundle $\LL$. We have the following well-known:

\begin{lemma}\label{17}   
$E\in SU_{X}(r,e)$ is not a base point for $|\LL^{p}|$ if there exists a vector 
bundle $F$ of rank $pr_{1}$ and degree $p(r_{1}(g-1)-e_{1})$ such that $h^{0}(E\otimes F)=0$.
\end{lemma}
It is easy to see that such an $F$ must necessarily be semistable (cf. \cite{Le Potier} 
(2.5)). It is also a simple consequence of the existence of Jordan-H\"older filtrations 
that one has to check the condition in the above lemma only for $E$ stable. We sketch 
the proof for convenience:

\begin{lemma}\label{18}
If for any stable bundle $V$ of rank $r^{\prime}\leq r$ and slope $e/r$ there exists 
$F\in U_{X}(pr_{1},p(r_{1}(g-1)-e_{1}))$ such that 
$h^{0}(V\otimes F)=0$, then the same is true for every $E\in SU_{X}(r,e)$. 
\end{lemma}
\begin{proof}
Assume that $E$ is strictly semistable. Then it has a Jordan-H\"older filtration:
$$0=E_{0}\subset E_{1}\subset\ldots\subset E_{n}=E$$
such that $E_{i}/E_{i-1}$ are stable for $i\in \{1,\ldots,n\}$ and $\mu (E_{i}/E_{i-1})=\frac{e}{r}$.
By assumption there exist $F_{i}\in U_{X}(pr_{1},p(r_{1}(g-1)-e_{1}))$ such that $h^{0}(E_{i}/E_{i-1}
\otimes F_{i})=0$ and so if we denote
$$\Theta_{E_{i}/E_{i-1}}:=\{F|h^{0}(E_{i}/E_{i-1}\otimes F)\neq 0\}\subset U_{X}(pr_{1},p(r_{1}(g-1)-e_{1})),$$
this is a proper subset for every $i$. It is clear that any 
$$F\in U_{X}(pr_{1},p(r_{1}(g-1)-e_{1}))- \underset{i=1}{\overset{n}{\bigcup}}\Theta_{E_{i}/E_{i-1}}$$
satisfies $h^{0}(E\otimes F)=0$.
\end{proof}
We also record a simple lemma which will be useful in \S4. It is
most certainly well known, but we sketch the proof for
convenience (cf. also \cite{Teixidor} 1.1).

\begin{lemma}\label{41}
Consider a sheaf extension:
$$0\longrightarrow F\longrightarrow E\longrightarrow
G\longrightarrow 0.$$
If $E$ is stable, then $h^{0}(G^{*}\otimes F)=0$.
\end{lemma}
\begin{proof}
Assuming the contrary, there is a nonzero morphism $G\rightarrow
F$. Composing this with the maps $E\rightarrow G$ to the left and 
$F\rightarrow E$ to the right, we obtain a nontrivial
endomorphism of $E$, which contradicts the stability assumption. 
\end{proof} 
As a final remark, note that we are always slightly abusing the notation by using vector bundles 
instead of $S$-equivalence classes. This is harmless, since it is easily seen that it is enough to 
check the assertions for any representative of the equivalence class.

\section{\textbf{An upper bound on the dimension of Hilbert schemes}}

The goal of this paragraph is to prove a result (see Theorem \ref{1} below) 
giving an upper bound on the 
dimension of the Hilbert schemes of coherent quotients of fixed rank and
degree of a given vector bundle, optimal at least in the case corresponding to line subbundles. 
For the general theory of Hilbert schemes the reader 
can consult for example \cite{Le Potier} \S4.  

Concretely, fix a vector bundle $E$ of rank $r$ and degree $e$ on $X$ and denote by 
${\rm Quot}_{r-k,e-d}(E)$ the Hilbert scheme of coherent quotients of $E$ of rank $r-k$ and 
degree $e-d$. 
We can (and will) identify ${\rm Quot}_{r-k,e-d}(E)$ to the 
set of subsheaves of $E$ of rank $k$ and degree $d$. Consider also:
$$d_{k}:={\rm max}\{\deg (F)|~F\subset E, {\rm rk}(F)=k\}$$
and
$$M_{k}(E)=\{F|~F\subset E, {\rm rk}(F)=k,\deg (F)=d_{k}\}$$
the set of maximal subbundles of rank $k$. Clearly any $F\in M_{k}(E)$ has to be a vector 
subbundle of $E$. 
Note that the number $e-d_{k}$ is exactly the minimal degree of a quotient bundle of $E$ 
of rank $r-k$.
By \cite{Mukai} \S2 we have the following basic result:

\begin{proposition}\label{2}
The following hold and are equivalent:
\newline
(i) {\rm dim~}$M_{k}(E)\leq k(r-k)$
\newline
(ii) for any $x\in X$ and any $W\subset E(x)$ $k$-dimensional subspace of the fiber of 
$E$ at $x$, there are at most finitely many $F\in M_{k}(E)$ such that $F(x)=W$.
\end{proposition} 

Part $(i)$ above thus gives an upper bound for the dimension of the  
Hilbert scheme in the case $d=d_{k}$. The next result is a generalization in the case of 
arbitrary degree $d$, which turns out to give an optimal result (see Example \ref{40} below).

\begin{theorem}\label{1}
With the notation above, we have:
$${\rm dim~}{\rm Quot}_{r-k,e-d}(E)\leq k(r-k) + (d_{k}-d)k(r-k+1).$$
\end{theorem}

\begin{remark}
To avoid any confusion, we emphasize here that the notation is slightly different from that 
used in the introduction, in the sense that we are replacing $k$ by $r-k$, $d$ by $e-d$ 
and $f_{k}$ by $e-d_{k}$. 
This is done for consistency in rewriting everything in terms of subbundles, but note that
the statement is exactly the same.
\end{remark}

The proof will proceed by induction on 
the difference $d_{k}-d$. In order to perform this induction we have to use a special case 
of the notion 
of \emph{elementary transformation} along a zero-dimensional subscheme of arbitrary 
length. We call this construction simply elementary transformation since there is no danger 
of confusion. 

\begin{definition}
Let $\tau$ be a zero-dimensional subscheme of $X$ supported on the points ${P_{1},\ldots,P_{s}}$. An 
\emph{elementary transformation} of $E$ along $\tau$ is a vector bundle $E^{\prime}$
defined by a sequence of the form:
$$0\longrightarrow E^{\prime}\longrightarrow E\overset{\phi}{\longrightarrow}
\tau\longrightarrow 0.$$
where the morphism $\phi$ is determined by giving surjective maps $E_{P_{i}}\overset{\phi_{i}}
{\rightarrow}\CC_{P_{i}}^{a_{i}}$ induced by specifying $a_{i}$ distinct hyperplanes in $E(P_{i})$
(whose intersection is the kernel of $\phi_{i}$)
, $\forall i\in \{1,\ldots,s\}$. We call $m=a_{1}+\ldots +a_{s}$ the \emph{length} of $\tau$
and $a_{i}$ the \emph{weight} of $P_{i}$.
\end{definition}
Let us briefly remark that this is not the most general definition, since we are imposing 
a condition on the choice of hyperplanes. We prefer to work with this notion because it is sufficient 
for our purposes and allows us to avoid some technicalities. Note though that the space parametrizing
these transformations is not compact. One could equally well work with the general definition,
when the hyperplanes could come together, and obtain a compact parameter space, which can be shown to 
be irreducible (it is basically a Hilbert scheme of rank zero quotients of fixed length). 

In fact it is an immediate observation that the elementary transformations of $E$ of length $m$, 
in the sense of the definition above, are parametrized by 
$Y:=(\PP E)_{m}- \Delta$, where $(\PP E)_{m}$ is the $m$-th
symmetric product of the projective bundle $\PP E$ and $\Delta$ is the union of all its 
diagonals. There is an obvious forgetful map 
$$\pi:Y\longrightarrow X_{m},$$
where $X_{m}$ is the $m$-th symmetric product of the curve $X$. We will denote by $Y_{m}\subset 
(\PP E)_{m}- \Delta$ the open subset $(\PP E)_{m}- \pi^{-1}(\delta)$, where 
$\delta$ is the union of the diagonals in $X_{m}$. Its points correspond to the elementary
transformations of length $m$ supported at $m$ distinct points of $X$.

\begin{definition}
Let $V$ be a subbundle of $E$. An elementary transformation 
$$0\longrightarrow E^{\prime}\longrightarrow E\overset{\phi}{\longrightarrow}
\tau\longrightarrow 0$$
is said to \emph{preserve} $V$ if the inclusion $V\subset E$ factors through the inclusion 
$E^{\prime}\subset E$.
\end{definition}

\begin{lemma}\label{42}
If $E^{\prime}$ is determined by the hyperplanes $H_{i}^{1},\ldots, H_{i}^{a_{i}}\subset E(P_{i})$
for $i\in \{1,\ldots,s\}$ and $V_{i}:=\overset{a_{i}}{\underset{j=1}{\bigcap}}H_{i}^{j}$, then $V$ 
is preserved by $E^{\prime}$ if and only if $V(P)\subset V_{i}, ~\forall i\in \{1,\ldots,s\}$.  
\end{lemma}
\begin{proof}
We have an induced diagram 
$$\xymatrix{
& & V \ar@{^{(}->}[d] \ar[dr]^{\alpha} & & \\
0 \ar[r] & E^{\prime} \ar[r] & E \ar[r]^{\phi} & \tau \ar[r] & 0 }$$
where $\alpha$ is the composition of $\phi$ with the inclusion of $V$ in $E$. It is clear that 
$E^{\prime}$ preserves $V$ if and only if $\alpha$ is identically zero. The lemma follows then 
easily from the definitions.
\end{proof}
In general it is important to know the dimension of the set of elementary transformations of a 
certain type preserving a given subbundle. This is given by the following simple lemma:

\begin{lemma}\label{3}
Let $V\subset E$ be a subbundle of rank $k$. Consider the set of elementary 
transformations of $E$ along a zero dimensional subscheme of length $m$ belonging to an 
irreducible subvariety $W$ of $X_{m}- \delta$ that preserve $V$:
$$\mathcal{Z}_{V}:=\{E^{\prime}|~V\subset E^{\prime}\}\subset \pi^{-1}(W),$$
where $\pi:Y_{m}\rightarrow X_{m}$ is the natural projection.
Then $\mathcal{Z}_{V}$ is irreducible of dimension $m(r-k-1)+{\rm dim~}W$.
\end{lemma}
\begin{proof}
An elementary transformation at $m$ points $x_{1},\ldots, x_{m}$ is given 
by a choice of hyperplanes $H_{i}\subset E(x_{i})$ for each $i$. By the previous lemma, 
such a transformation 
preserves $V$ if and only if $V(x_{i})\subset H_{i}$ for all $i$. 
We have a natural diagram
$$\xymatrix{
\mathcal{Z}_{V} \ar[dr]_{p} \ar@{^{(}->}[r]^{i} & \pi^{-1}(W) \ar[d]^{\pi} \\
& W }$$
where $\pi$ is the restriction to $\pi^{-1}(W)$ and $p$ is the 
composition with the inclusion of $\mathcal{Z}_{V}$ in $\pi^{-1}(W)$. For $D=x_{1}+\ldots +
x_{m}\in W$, we have:
$$p^{-1}(D)\cong \{(H_{1},\ldots ,H_{m})|~V(x_{i})\subset H_{i}\subset E(x_{i}),\forall 
i=1,\ldots ,m\}$$
$$\cong \PP^{r-k-1}\times\ldots\times \PP^{r-k-1}$$
where the product is taken $m$ times. So $p^{-1}(D)$ is irreducible of dimension $m(r-k-1)$
and this gives that $\mathcal{Z}_{V}$ is irreducible of dimension $m(r-k-1)+{\rm dim~}W$.
\end{proof}
The following proposition will be the key step in running the inductive argument.
It computes ``how fast'' we can eliminate all the maximal subbundles of $E$ while preserving 
a fixed nonmaximal subbundle.

\begin{proposition}\label{4}
Let $V\subset E$ be a subbundle of rank $k$ and degree $d$, not maximal. Then if $m\geq 
r-k+1$, there exists an elementary transformation of length $m$ 
$$0\longrightarrow E^{\prime}\longrightarrow E\longrightarrow \tau\longrightarrow 0$$
such that $V\subset E^{\prime}$, but $F\not\subset E^{\prime}$ for any $F\in M_{k}(E)$.
In other words $E^{\prime}$ preserves $V$, but does not preserve any maximal $F$. 
If we fix a point $P\in X$, then $\tau$ can be chosen to have weight $m-1$ at $P$ and
weight $1$ at a generic point $Q\in X$. 
\end{proposition} 
\begin{proof}
Fix a point $P\in X$. We can consider an elementary transformation of $E$ of length 
$r-k$, supported only at $P$:
$$0\longrightarrow E^{\prime}\longrightarrow E\longrightarrow \tau\longrightarrow 0,$$ 
such that ${\rm Im}(E^{\prime}(P)\rightarrow E(P))=V(P)$. Then as
in \ref{42}, the only maximal 
subbundles $F$ that are preserved by this transformation are exactly those such that 
$F(P)=V(P)$. By Proposition \ref{2} this implies that only at most a finite number of $F$'s 
can be preserved. 

If none of the maximal subbundles actually survive in $E^{\prime}$, then any further 
transformation at one point would do. Otherwise clearly for a generic $Q\in X$ we have 
$F(Q)\neq V(Q)$ for all the $F$'s that are preserved and so we can choose a hyperplane 
$V(Q)\subset H\subset E(Q)$ such that $F(Q)\not\subset H$ for any such $F$. The 
elementary transformation of $E^{\prime}$ at $Q$ corresponding to this hyperplane 
satisfies then the required property.  
\end{proof}

\begin{remark}\label{21}
(1) It can definitely happen that all the maximal subbundles are killed by elementary 
transformations of length less than $r-k+1$ which preserve $V$. In any case, as it was 
already suggested in the proof above, by further elementary transforming we obviously 
don't change the property that we are interested in, so $r-k+1$ is a bound that works 
in all situations. 
\newline
(2) By Lemma \ref{3}, the set $\mathcal{Z}_{V}$ of all elementary transformations of length $m$ preserving $V$ 
is irreducible of dimension $m(r-k)$. On the other hand the condition of preserving at 
least one maximal subbundle is closed, so once the lemma above is true for one elementary 
transformation, it applies for an open subset of $\mathcal{Z}_{V}$.
\end{remark}
Finally we have all the ingredients necessary to prove the theorem.
To simplify the formulations, it is convenient to introduce the following ad-hoc definition:

\begin{definition}
An irreducible component $\mathcal{Q}\subset {\rm Quot}_{r-k,e-d}(E)$ is called
\emph{non-special} if its generic point corresponds to a locally free quotient of $E$ 
and \emph{special} if all its points correspond to non-locally free quotients.
For any $\mathcal{Q}$, denote by $\mathcal{Q}_{0}$ the open subset parametrizing 
locally free quotients and consider ${\rm Quot}_{r-k,e-d}^{0}(E):=\underset{\mathcal{Q}}
{\bigcup}\mathcal{Q}_{0}$. 
\end{definition}

\begin{proof}(of \ref{1})
Denote by $\mathcal{Q}$ an irreducible component of ${\rm Quot}_{r-k,e-d}(E)$ (recall that 
we are thinking now of this Hilbert scheme as parametrizing subsheaves of rank $k$ and 
degree $d$). The first step is to observe that it is enough to prove the statement 
when $\mathcal{Q}$ is non-special. To see this, note that every nonsaturated subsheaf $F\subset E$
determines a diagram:
$$\xymatrix{
& & 0 \ar[d] & 0 \ar[d] \\
0 \ar[r] & F \ar[r] \ar[d]_{\cong} & F^{\prime} \ar[r] \ar[d] & \tau \ar[r]
\ar[d] & 0 \\
0 \ar[r] & F \ar[r] & E \ar[r] \ar[d] & G^{\prime} \ar[r] \ar[d] & 0 \\
& & G \ar[d] \ar[r]^{\cong} & G \ar[d] \\
& & 0 & 0 } $$
where $F^{\prime}$ is the saturation of $F$, $G$ is a quotient vector bundle and $\tau$, the torsion 
subsheaf of $G^{\prime}$, is a
nontrivial zero-dimensional subscheme, say of length $a$. We can stratify the set of all such $F$'s 
according to the value of the parameter $a$, which obviously runs over a finite set. If we denote 
by $\{F\}_{a}$ the subset corresponding to a fixed $a$, this gives then:
$${\rm dim~}\{F\}_{a}\leq {\rm dim~}{\rm Quot}_{r-k,e-d-a}^{0}(E)+ka.$$
The right hand side is clearly less than $k(r-k)+(d_{k}-d)k(r-k+1)$ if we assume that the 
statement of the theorem holds for ${\rm Quot}_{r-k,e-d-a}^{0}(E)$. 

Let us then restrict to the case when $\mathcal{Q}$ is a non-special component.
The proof goes by induction on $d_{k}-d$. If $d_{k}=d$, the statement is exactly the 
content of \ref{2}. Assume now that $d_{k}>d$ and that the statement holds for all the
pairs where this difference is smaller.
Recall that $\mathcal{Q}_{0}\subset \mathcal{Q}$ denotes the open subset corresponding to vector 
subbundles
and fix $V\in \mathcal{Q}_{0}$. Then by Proposition \ref{4}, there exists an 
elementary transformation  
$$0\longrightarrow E^{\prime}\longrightarrow E\longrightarrow \tau\longrightarrow 0$$
of length $r-k+1$, such that $V\subset E^{\prime}$, but $F\not\subset E^{\prime}$ for any $F\in M_{k}(E)$.
Then necessarily $d_{k}(E^{\prime})<d_{k}(E)=d_{k}$ (consider the saturation in $E$ of 
a maximal subbundle of $E^{\prime}$) and so $d_{k}(E^{\prime})-d<d_{k}-d$. This means 
that we can apply the inductive hypothesis for any non-special component of the set of subsheaves
of rank $k$ and degree $d$ of $E^{\prime}$.
To this end, consider the correspondence:

\xymatrix{ \\
& \mathcal{W} \ar[dl]_{p_{1}} \ar[dr]^{p_{2}} & = & \{(V,E^{\prime})|~V\subset E^{\prime},
F\not\subset E^{\prime},\forall F\in M_{k}(E)\}\subset \mathcal{Q}_{0}\times Y_{r-k+1} \\
\mathcal{Q}_{0} & & Y_{r-k+1}. }
By Lemma \ref{3} and Remark \ref{21}(b), for any 
$V\in \mathcal{Q}_{0}$, the fiber 
$p_{1}^{-1}(V)$ is a (quasi-projective irreducible) variety of dimension $(r-k+1)(r-k)$ and so:
\begin{equation}
{\rm dim~} \mathcal{W}={\rm dim~} \mathcal{Q}_{0}+(r-k+1)(r-k).
\end{equation}
On the other hand, for $E^{\prime}\in {\rm Im}(p_{2})$, the inductive hypothesis implies that 
$${\rm dim~} p_{2}^{-1}(E^{\prime})\leq k(r-k)+(d_{k}(E^{\prime})-d)k(r-k+1)$$
$$\leq k(r-k)+(d_{k}-d-1)k(r-k+1)$$
and since ${\rm dim~}Y_{r-k+1}=r(r-k+1)$ we have:
\begin{equation}
{\rm dim~}\mathcal{W}\leq r(r-k+1)+k(r-k)+(d_{k}-d-1)k(r-k+1).
\end{equation}
Combining (1) and (2) we get:
$${\rm dim~}\mathcal{Q}_{0}\leq k(r-k)+(d_{k}-d)k(r-k+1)$$
and of course the same holds for $\mathcal{Q}=\overline{\mathcal{Q}_{0}}$. This completes
the proof. 
\end{proof}
The formulation and the proof of the theorem give rise to a few natural questions and we address them in the 
following examples.

\begin{example}
It is easy to construct special components of Hilbert schemes. For example consider for any 
$X$ the Hilbert scheme of quotients of $\OO_{X}^{\oplus 2}$ of rank $1$ and degree $1$. There certainly 
exist such quotients which have torsion, like
$$\OO_{X}^{\oplus 2}\longrightarrow \OO_{X}\oplus \OO_{P}\longrightarrow 0,$$
where $P$ is any point of $X$, but for obvious cohomological reasons there can be no sequence of the form 
$$0\longrightarrow L^{-1}\longrightarrow \OO_{X}^{\oplus 2}\longrightarrow L\longrightarrow 0$$
with ${\rm deg}(L)=1$. So in this case there are actually no non-special components.
\end{example} 

\begin{example}
Going one step further, there may exist special components whose dimension is greater 
than that of any of the non-special ones. Note though that the proof shows that in this 
case that the bound cannot be optimal. To see an example, consider quotients of $\OO_{X}^{\oplus 2}$
of rank $1$ and degree $1\leq d\leq g-2$ on a nonhyperelliptic curve $X$. Any such locally free
quotient $L$ gives a sequence:
$$0\longrightarrow L^{-1}\longrightarrow \OO_{X}^{\oplus 2}\longrightarrow L\longrightarrow 0$$
and so the dimension of any component of the Hilbert scheme containing it is bounded above 
by $h^{0}(L^{\otimes 2})$. Now Clifford's theorem says that $h^{0}(L^{\otimes 2})\leq d+1$, but 
our choices make the equality case impossible, so in fact $h^{0}(L^{\otimes 2})\leq d$.

On the other hand consider an effective divisor $D$ of degree $d$. Then a point in the same 
Hilbert scheme is determined by a natural sequence:
$$0\longrightarrow \OO_{X}(-D)\longrightarrow \OO_{X}^{\oplus 2}\longrightarrow 
\OO_{X}\oplus\OO_{D}\longrightarrow 0$$ 
and it is not hard to see that the dimension of the Hilbert scheme at this point is equal 
to $d+1$ (essentially $d$ parameters come from $D$ and one from the sections of $\OO_{X}^{\oplus 2}$).
This gives then a special component whose dimension is greater than that of any non-special one.
\end{example}

\begin{example}\label{40}
More significantly, the bound given in the theorem is optimal. Consider for this 
a line bundle $L$ of degree $4$ on a curve $X$ of genus $2$ and a generic extension:
$$0\longrightarrow \OO_{X}\longrightarrow E\longrightarrow L\longrightarrow 0.$$
By standard arguments one can see that such an extension must be stable. Since $\mu(E)=2$, by the 
classical theorem of Nagata \cite{Nagata} we get that $d_{1}(E):=\underset{M\subset E}{{\rm max}} 
{\rm~ deg}(M)=1$
and so for the sequence above $d_{1}-d =1$. The theorem then tells us that the dimension of any component 
of the Hilbert scheme containing the given quotient is bounded above by $3$. But on the other hand 
$h^{1}(L)=0$, so this gives a smooth point and the dimension of the component is $h^{0}(L)$, which 
by Riemann-Roch is exactly $3$.

This example turns out to be a special case of a general pattern, as suggested by M. Teixidor. In fact 
in \cite{Teixidor} it is shown that whenever $E$ is a generic stable bundle, the invariant $d_{k}$
is the largest integer $d$ that makes the expression $ke-rd-k(r-k)(g-1)$ nonnegative (cf. also \cite{Brambila}).
Also the dimension of the Hilbert scheme can be computed exactly in this case (see \cite{Teixidor} 0.2),
and for example under the numerical assumptions above it is precisely equal to $3$. Thus in fact for 
every generic stable bundle of rank $2$ and degree $4$ on a curve of genus $2$, we have equality in the 
theorem. Much more generally, it can be seen analogously that for any $r$ and $g$ equality is satisfied
for a generic stable bundle as long as $d_{1}$ satisfies a certain numerical condition.
\end{example}

The proof of the theorem given above can be slightly modified towards a more 
natural and compact form. We chose to follow
the longer approach because it emphasizes very clearly what is the phenomenon involved, but below we would
also like to briefly sketch this parallel argument, which grew out of conversations with I. Coand\u a.
 
We will use the same notations as above. 
There exists a natural specialization map:
$$X\times {\rm Quot}_{r-k,e-d}^{0}(E)\longrightarrow {\bf{G}}_{r-k}(E),$$
where ${\bf{G}}_{r-k}(E)$ is the Grassmann bundle of $r-k$ dimensional quotients of the fibers of $E$. Of course
in the case $d=d_{k}$, ${\rm Quot}_{r-k,e-d}^{0}(E)$ is compact and the morphism above is finite. Fix now
$P\in X$ and $w\in {\bf{G}}_{r-k}(E(P))$ a point corresponding to a quotient $E(P)\rightarrow W
\rightarrow 0$. The choice of $P$ determines a map:
$$\phi: {\rm Quot}_{r-k,e-d}^{0}(E)\longrightarrow {\bf{G}}_{r-k}(E(P))$$
and we would like to bound the dimension of $\phi^{-1}(w)$.
There is a natural induced sequence:
$$0\longrightarrow F\longrightarrow E\longrightarrow W\otimes {\bf{C}}_{P}\longrightarrow 0$$
and it is not hard to see that $\phi^{-1}(w)$ embeds as an open subset in  
${\rm Quot}_{r-k,e-d-k}^{0}(F)$. Every locally free
quotient of $F$ has degree $\geq e-d_{k}-k$, and there are at most a finite number of quotients having
precisely this degree (they come exactly from the minimal degree quotients of $E$ having fixed fiber
$W$ at $P$). Let $G_{1},\ldots , G_{m}$ be these quotients, sitting in exact sequences:
$$0\longrightarrow F_{i}\longrightarrow F\longrightarrow G_{i}\longrightarrow 0.$$
The variety $Y:= \PP F - \underset{i=1}{\overset{m}{\bigcup}}\PP G_{i}$ parametrizes then the one-point
elementary transformations of $F$ that do not preserve any of the $F_{i}$'s. Consider the natural incidence
$$\mathcal{Z}\subset {\rm Quot}_{r-k,e-d-k}^{0}(F)\times Y$$
parametrizing the pairs $(F\rightarrow Q\rightarrow 0,F^{\prime})$, where $F^{\prime}$ is an elementary
transformation in $Y$ and $Q$ is not preserved as a quotient of $F^{\prime}$ (in other words the corresponding
kernel is preserved). The fiber of $\mathcal{Z}$ over $F\rightarrow Q\rightarrow 0$ is isomorphic to
$\PP Q\cap Y$ and so
$${\rm dim}~\mathcal{Z}={\rm dim}{\rm~ Quot}_{r-k,e-d-k}^{0}(F)+r-k.$$
On the other hand the fiber of $\mathcal{Z}$ over $F^{\prime}\in Y$ is ${\rm Quot}_{r-k,e-d-k-1}^{0}(F^{\prime})$
Now for $F^{\prime}$ the minimal degree of a quotient of rank $r-k$ is smaller, hence inductively as before:
$${\rm dim}{\rm~ Quot}_{r-k,e-d-k-1}^{0}(F^{\prime})\leq k(r-k)+(d_{k}-d-1)k(r-k+1).$$
This immediately implies that
$${\rm dim}{\rm ~Quot}_{r-k,e-d-k}^{0}(F)\leq k(r-k)+(d_{k}-d-1)k(r-k+1)+k.$$
As this consequently holds for every fiber of the map $\phi$, we conclude that
$${\rm dim}{\rm ~Quot}_{r-k,e-d}^{0}(E)\leq k(r-k)+(d_{k}-d)k(r-k+1),$$
which finishes the proof.

\section{\textbf{Warm up for effective base point freeness: the case of $SU_{X}(2)$}}
 
In this section we give a very simple proof of a theorem which first appeared in
\cite{Raynaud}
(see also \cite{Hein}). It completely takes care of the case of $SU_{X}(2)$.
Although the specific technique (based on the Clifford theorem for line bundles) is
different from the methods that will be used in Section 4 to prove the main result, the
general computational idea already appears here, in a particularly transparent form.
This is the reason for including the proof.
\begin{theorem}\label{6}
The linear system $|\LL|$ on $SU_{X}(2)$ has no base points.
\end{theorem}
\begin{proof}
Recall from \ref{17} and \ref{18} that
the statement of the theorem is equivalent to the following fact: for any stable bundle
$E\in SU_{X}(2)$, there exists a line bundle $L\in {\rm Pic}^{g-1}(X)$ such that   
$H^{0}(E\otimes L)=0$. This is certainly an open condition and it is sufficient to
prove that the algebraic set
$$\{L\in {\rm Pic}^{g-1}(X)|~H^{0}(E\otimes L)\neq 0\}\subset {\rm Pic}^{g-1}(X)$$
has dimension strictly less than $g$.
 
A nonzero map $E^{*}\rightarrow L$ comes together with a diagram of the form:
$$\xymatrix{ 
E^{*} \ar[r] \ar[dr] & M  \ar[r] \ar@{^{(}->}[d] & 0 \\
 & L }$$
where $M$ is just the image in $L$. Then we have $M=L(-D)$ for some effective      
divisor $D$. Since $E$ is stable, the degree of $M$ can vary from $1$ to $g-1$ and
we want to count all these cases separately. So for $m=1,\ldots , g-1$, consider
the following algebraic subsets of ${\rm Pic}^{g-1}(X)$:
$$A_{m}:=\{L\in {\rm Pic}^{g-1}(X)|~\exists~ 0\neq\phi:E^{*}\rightarrow L {\rm ~with~}
M={\rm Im}(\phi),{\rm deg}(M)=m\}.$$
The claim is that ${\rm dim~}A_{m}\leq g-1$ for all such $m$. Then of course 
$$A_{1}\cup\ldots\cup A_{g-1}\subsetneq {\rm Pic}^{g-1}(X)$$
and any $L$ outside this union satisfies our requirement.
To prove the claim, denote by ${\rm Quot}_{1,m}(E)$ the Hilbert scheme of coherent
quotients of
$E$ of rank $1$ and degree $m$. The set of line bundle quotients $E^{*}\rightarrow
M\rightarrow 0$ of degree $m$ is a subset of ${\rm Quot}_{1,m}(E)$. On the
other hand every $L\in A_{m}$ can be written as $L=M(D)$, with $M$ as above and
$D$ effective of degree $g-1-m$. This gives the obvious bound:
$${\rm dim~}A_{m}\leq {\rm dim~}{\rm Quot}_{1,m}(E)+g-1-m.$$
To bound the dimension of the Hilbert scheme in question, fix an $M$ as before
and consider the exact sequence that it determines:
$$0\longrightarrow M^{*}\longrightarrow E\longrightarrow M\longrightarrow 0.$$
Note that the kernel is isomorphic to $M^{*}$ since $E$ has trivial determinant. 
Now we use the well known fact from deformation theory that ${\rm dim~}{\rm Quot}_{1,m}(E)
\leq h^{0}(M^{\otimes 2})$.
To estimate $h^{0}(M^{\otimes 2})$, one uses all the information provided 
by Clifford's theorem. The initial bound that it gives is $h^{0}(M^{\otimes 2})
\leq m+1$ (note that ${\rm deg}(M)\leq g-1$). 
If actually $h^{0}(M^{\otimes 2})\leq m$, then we immediately
get ${\rm dim~}A_{m}\leq g-1$ as required. On the other hand if
$h^{0}(M^{\otimes 2})=m+1$, by the equality case in Clifford's theorem (see e.g. \cite{ACGH} III, \S1) 
one of the following must hold:
$M^{\otimes 2}\cong
\OO_{X}$ or $M^{\otimes 2}\cong \omega_{X}$ or $X$ is hyperelliptic and
$M^{\otimes 2}\cong m\cdot g_{2}^{1}$. The first case is impossible since
${\rm deg}(M)>0$. In the second case $M$ is a theta characteristic and
we are done either by the fact that these are a finite number or by other
overlapping cases. The third case can also happen only for a finite number of
$M$'s and if we're not in any of the other cases then of course ${\rm dim~}A_{m}\leq
g-1-m<g-1$. This concludes the proof of the theorem.
\end{proof}
\begin{remark}
Note that the key point in the proof above is the ability to give a convenient
upper bound on the dimension of certain Hilbert schemes. This will essentially
be the main ingredient in the general result proved in Section 4, and the needed
estimate was provided in the previous section.
\end{remark}

\section{\textbf{Base point freeness for pluritheta linear series on $SU_{X}(r,e)$}}

Using the dimension bound given in Section 2, we are now able to prove the main
result of this paper, namely an effective base point freeness bound for 
pluritheta linear series on $SU_{X}(r,e)$. The proof is computational in 
nature and the roots of the main technique involved have already appeared in \ref{6}.
Let $r\geq 2$ and $e$ be arbitrary integers and let $h={\rm
gcd}(r,e),~r=r_{1}h$ and $e=e_{1}h$. For the statement it is
convenient to introduce another invariant of the moduli space. If
$E\in SU_{X}(r,e)$ and $1\leq k\leq r-1$, define
$s_{k}(E):=ke-rd_{k}$, where $d_{k}$ is the maximum degree of a 
subbundle of $E$ of rank $k$ (cf.\cite{Lange2}). Note that if $E$
is stable one has $s_{k}(E)\geq h$ and we can further define 
$s_{k}=s_{k}(r,e):=\underset{E~ {\rm stable}}{{\rm min}}s_{k}(E)$
and $s=s(r,e):=\underset{1\leq k\leq r-1}{{\rm min}}s_{k}$. 
Clearly $s\geq h$ and it is also an immediate
observation that $s(r,e)=s(r,-e)$.

\begin{theorem}\label{10}
The linear series $|\LL^{p}|$ on $SU_{X}(r,e)$ is base point free
for 
$$p\geq {\rm  max}\{
\frac{(r+1)^{2}}{4r} h,\frac{r^{2}}{4s}h\}.$$ 
\end{theorem}
\begin{remark}
Note that the bound given in the theorem is always either a
quadratic or a linear function in the rank $r$.
It should also be said right away that although this bound works uniformly, in 
almost any particular situation one can do a little better. Unfortunately there doesn't seem 
to be a better uniform way to express it, but we will comment more on this at the end of 
the section (cf. Remark \ref{44}). 
\end{remark}

\begin{proof}(of \ref{10})
Let us denote for simplicity $M:={\rm  max}\{
\frac{(r+1)^{2}}{4r} h,\frac{r^{2}}{4s}h\}$.
Since the problem depends only on the residues of $e$ modulo $r$, there is no loss of 
generality in looking only at $SU_{X}(r,-e)$ with $0\leq e\leq r-1$.
The statement of the theorem is implied by the following assertion, as described in \ref{17} and 
\ref{18}:
$$ For~ any~stable~ E\in SU_{X}(r,-e)~ and~ any~ p\geq M,~there~is~an $$
$$F\in U_{X}(pr_{1},p(r_{1}(g-1)+e_{1}))~such~that~h^{0}(E\otimes F)=0.$$ 

Fix a stable bundle $E\in SU_{X}(r,-e)$. 
Note that only in this proof, as opposed to the rest of the paper, $e$ in fact denotes 
the degree of $E^{*}$, and not that of $E$.
If for some $F\in U_{X}(pr_{1},p(r_{1}(g-1)+e_{1}))$ there is a 
nonzero map  $E^{*}\overset{\phi}{\rightarrow} F$, then this comes together with a 
diagram of the form: 
$$\xymatrix{
E^{*} \ar[r] \ar[dr]_{\phi} &  V \ar[r] \ar@{^{(}->}[d] & 0 \\
 & F }$$
where the vector bundle $V$ is the image of $\phi$. The idea is essentially to count 
all such diagrams assuming that the rank and degree of $V$ are fixed and see that the 
$F$'s involved in at least one of them 
cover only a proper subset of the whole
moduli space. Denote as before by ${\rm Quot}_{k,d}(E^{*})$ the Hilbert scheme of quotients 
of $E^{*}$ of rank $k$ and degree $d$ and for any $1\leq k\leq r$ and any $d$ in the suitable 
range (given by the stability of $E$ and $F$) consider its subset:
$$A_{k,d}:=\{V\in {\rm Quot}_{k,d}(E^{*})|~\exists F\in U_{X}(pr_{1},p(r_{1}(g-1)+e_{1})),
~\exists 0\neq\phi:E^{*}\rightarrow F 
{\rm~ with~} V={\rm Im}(\phi)\}.$$
The theorem on Hilbert schemes stated in the introduction then gives us the dimension estimate:
\begin{equation}
{\rm dim~}A_{k,d}\leq k(r-k)+(d-f_{k})(k+1)(r-k),
\end{equation}
where $f_{k}=f_{k}(E^{*})$ is the minimum possible degree of a quotient bundle of $E^{*}$ of rank $k$
(which is the same as $-d_{k}$).
Define now the following subsets of $U_{X}(pr_{1},p(r_{1}(g-1)+e_{1}))$:
$$U_{k,d}:=\{F|~\exists V\in A_{k,d} {\rm ~with~} V\subset F\}\subset U_{X}(pr_{1},p(r_{1}(g-1)+e_{1})).$$
The elements of $U_{k,d}$ are all the $F$'s that appear in diagrams as above for 
fixed $k$ and $d$. The claim is that 
$${\rm dim~}U_{k,d}< (pr_{1})^{2}(g-1)+1,$$
which would imply that $U_{k,d}\subsetneq U_{X}(pr_{1},p(r_{1}(g-1)+e_{1}))$. Assuming that this is true,
and since $k$ and $d$ run over 
a finite set, any $F\in U_{X}(pr_{1},p(r_{1}(g-1)+e_{1}))- \underset{k,d}{\bigcup}U_{k,d}$ 
satisfies the desired property that $h^{0}(E\otimes F)=0$, which gives the statement 
of the theorem. It is easy to see, and in fact a particular case
of the computation below, that in the case $k=r$ (i.e. $V=E^{*}$)
$U_{k,d}$ has dimension exactly $(pr_{1})^{2}(g-1)$.   

Let us concentrate then on proving the claim above for $1\leq
k\leq r-1$. Note that the inclusions 
$V\subset F$ appearing in the definition of $U_{k,d}$ are valid in general 
only at the sheaf level. Any such inclusion determines an exact sequence:
\begin{equation}
0\longrightarrow V\longrightarrow F\longrightarrow G^{\prime}\longrightarrow 0,
\end{equation}
where $G^{\prime}=G\oplus \tau_{a}$, with $G$ locally free and  
$\tau_{a}$ a zero dimensional subscheme of length $a$.
We stratify $U_{k,d}$ by the subsets 
$$U_{k,d}^{a}:=\{F| ~F ~{\rm given~by~ an ~extension~ of ~type~(4)}\}\subset U_{k,d},$$
where $a$ runs over the obvious allowable finite set of integers. 
A simple computation shows that $G$ has rank $pr_{1}-k$ and degree $p(r_{1}(g-1)+e_{1})-d-a$. Denote 
by $T_{k,d}^{a}$ the set of all vector bundles $G$ that are quotients of 
some $F\in U_{X}(pr_{1},p(r_{1}(g-1)+e_{1}))$. These can be parametrized by a relative Hilbert scheme 
(see e.g. \cite{Le Potier2} \S8.6)
over (an \'etale cover of) $U_{X}(pr_{1},p(r_{1}(g-1)+e_{1}))$ 
and so they form a bounded family. We invoke a general result, proved in 
\cite{BGN} 4.1 and 4.2,
saying that the dimension of such a family is always at most what we get if we assume 
that the generic member is stable. Thus we get the bound:
$${\rm dim~}T_{k,d}^{a}\leq (pr_{1}-k)^{2}(g-1)+1.$$ 

Now we only have to compute the dimension of the family of all possible extensions 
of the form (4) when $V$ and $G$ are allowed to vary over $A_{k,d}$ and 
$T_{k,d}^{a}$ respectively and $\tau_{a}$ varies over the symmetric product 
$X_{a}$. Any such extension induces a diagram
$$\xymatrix{
& & 0 \ar[d] & 0 \ar[d] \\
0 \ar[r] & V \ar[r] \ar[d]_{\cong} & V^{\prime} \ar[r] \ar[d] & \tau_{a} \ar[r]
\ar[d] & 0 \\
0 \ar[r] & V \ar[r] & F \ar[r] \ar[d] & G^{\prime} \ar[r] \ar[d] & 0 \\
& & G \ar[d] \ar[r]^{\cong} & G \ar[d] \\
& & 0 & 0 }$$ 
If we denote by $A_{k,d}^{a}$ the set of isomorphism classes of vector bundles 
$V^{\prime}$ that are (inverse) elementary transformations of length $a$ of 
vector bundles in $A_{k,d}$, then we have the obvious:
$${\rm dim~}A_{k,d}^{a}\leq {\rm dim~}A_{k,d}+ka.$$
On the other hand any $F$ is obtained as an extension of a bundle in 
$T_{k,d}^{a}$ by a bundle in $A_{k,d}^{a}$. Denote by $\mathcal{U}
\subset A_{k,d}^{a}\times T_{k,d}^{a}$ 
the open subset consisting of pairs $(V^{\prime},G)$ such that there exists an 
extension
$$0\longrightarrow V^{\prime}\longrightarrow F\longrightarrow G\longrightarrow 0$$
with $F$ stable. Note that by Lemma \ref{41} for any such pair we have $h^{0}(G^{*}\otimes 
V^{\prime})=0$ and so by Riemann-Roch $h^{1}(G^{*}\otimes V^{\prime})$ is constant,
given by:
\begin{equation}
h^{1}(G^{*}\otimes V^{\prime})=2kpr_{1}(g-1)-k^{2}(g-1)+kpe_{1}-pr_{1}d-pr_{1}a.
\end{equation}
In this situation it is a well known result (see e.g. \cite{Ramanan} (2.4) or 
\cite{Lange} \S4) that there exists a universal space of extension classes 
$\PP (\mathcal{U})\rightarrow \mathcal{U}$ whose dimension is computed by the 
formula:
$${\rm dim~}\PP (\mathcal{U})={\rm dim~}A_{k,d}^{a}+{\rm dim~}T_{k,d}^{a}
+h^{1}(G^{*}\otimes V^{\prime})-1.$$
There is an obvious forgetful map:
$$\PP (\mathcal{U})\longrightarrow U_{X}(pr_{1},p(r_{1}(g-1)+e_{1}))$$
whose image is exactly $U_{k,d}^{a}$. Thus by putting together
all the inequalities above we obtain:
$${\rm dim~}U_{k,d}^{a}\leq {\rm dim~}A_{k,d}^{a}+{\rm dim~}T_{k,d}^{a}
+h^{1}(G^{*}\otimes V^{\prime})-1$$
$$\leq k(r-k)+(d-f_{k})(k+1)(r-k)+ka +(pr_{1}-k)^{2}(g-1)+1$$
$$+2kpr_{1}(g-1)-k^{2}(g-1)+kpe_{1}-pr_{1}d-pr_{1}a-1$$
$$\leq k(r-k)+(d-f_{k})(k+1)(r-k)+(pr_{1})^{2}(g-1)+kpe_{1}-pr_{1}d+ka-pr_{1}a$$
$$\leq k(r-k)+(d-f_{k})(k+1)(r-k)+(pr_{1})^{2}(g-1)+kpe_{1}-pr_{1}d,$$
where the last inequality is due to the obvious fact that $k\leq
r-1< pr_{1}$ if $p\geq M$. 
Since $a$ runs over a finite set, to conclude the proof of the claim it is enough 
to see that ${\rm dim~}U_{k,d}^{a}\leq (pr_{1})^{2}(g-1)$. By the inequality above this 
is true if 
$$p(r_{1}d-ke_{1})\geq k(r-k) + (d-f_{k})(k+1)(r-k),$$
or equivalently if 
$$p(rd-ke)\geq k(r-k)h + (d-f_{k})(k+1)(r-k)h$$
for any $k$ and $d$. 
This can be rewritten in the following more manageable form:
$$p(r(d-f_{k})+rf_{k}-ke)\geq k(r-k)h + (d-f_{k})(k+1)(r-k)h.$$
The first case to look at is $d=f_{k}$, when we should have $p(rf_{k}-ke)\geq k(r-k)h$ 
and this should hold for every $k$.
But clearly $rf_{k}-ke=s_{r-k}(E^{*})=s_{k}(E)$, defined above in terms of
maximal subbundles, and $h|s_{k}(E)$. Since $E$ is stable we then
have $s_{k}(E)\geq h$ for all $k$ and so $s\geq h$ as mentioned
before. Note though that in general one cannot do better
(cf. Remark \ref{61}). In any case, this says that the inequality
$p\geq \frac{r^{2}}{4s} h$ must be satisfied (which would certainly
hold if $p\geq \frac{r^{2}}{4}$).

When $d>f_{k}$, it is convenient to collect together all the terms containing $d-f_{k}$. The last 
inequality above then reads:
$$(d-f_{k})(pr-(k+1)(r-k)h)+ps_{k}(E)\geq k(r-k)h.$$
For $p$ as before it is then sufficient to have $pr\geq
(k+1)(r-k)h$, which again by simple optimization is satisfied for
$p\geq \frac{(r+1)^{2}}{4r} h$. Concluding, the desired
inequality holds as long as $p\geq M$. 
\end{proof}

The most important instances of this theorem are the cases
of vector bundles of degree $0$ (more generally $d\equiv 0({\rm
  mod}~r)$) and degree $1$ or $-1$ (more generally $d\equiv
\pm 1 ({\rm mod}r)$). In the second situation the
moduli space in question is smooth. It is somewhat surprising
that the results obtained in these cases have different orders of 
magnitude.

\begin{corollary}\label{30}
$|\LL^{p}|$ is base point free on $SU_{X}(r)$ for $p\geq \frac{(r+1)^{2}}{4}$.
\end{corollary}
\begin{proof}
This is clear since $h=r$.
\end{proof}

\begin{corollary}\label{31}
$|\LL^{p}|$ is base point free on $SU_{X}(r,1)$ and
$SU_{X}(r,-1)$ for $p\geq r-1$.
\end{corollary}
\begin{proof}
Note that by duality it
suffices to prove the claim for one of the moduli spaces, say
$SU_{X}(r,-1)$.
In this case $h=1$ and for any $E\in SU_{X}(r,-1)$,
$s_{r-k}(E^{*})=rf_{k}-k\geq r-k$. Following the proof of the
theorem we see thus that it suffices to have 
$$p\geq {\rm max}\{r-1,\frac{(r+1)^{2}}{4r}\},$$
which is equal to $r-1$ for $r\geq 3$. For $r=2$ one can slightly
improve the last inequality in the proof of the theorem (actually
this is true whenever $r$ is even) to see that $p=1$ already works. 
\end{proof}

\begin{remark}\label{44}
Corollary \ref{31} can be seen as a generalization of the 
well-known fact that $|\LL|$ is base point free 
on $SU_{X}(2,1)$. Also, as already noted in its proof, the
general bound obtained in the theorem can be slightly improved in 
each particular case, due to the fact that the two optimization 
problems do not simultaneously have integral solutions. Thus for
example if $r$ is even, the proof of the theorem actually gives that $|\LL^{p}|$
is base point free on $SU_{X}(r)$ for $p\geq r(r+2)/4$.  
\end{remark}

\begin{remark}\label{61}
As already noted, in any given numerical situation the bound
given by the theorem is either linear or quadratic in the rank
$r$. One may thus hope that at least in the case $h=1$ (i.e. $r$
and $e$ coprime), a closer study of the number $s(r,e)$ might
always produce by this method a linear bound. Examples show
though that this is not the case: one can take $r=4l$, $e=2l-1$, 
$k=2l+1$ and $f_{k}=l$ (this works for special vector bundles)
for a positive integer $l$, which implies $s=s_{k}=1$. 
\end{remark}

\section{\textbf{Effective base point freeness on $U_{X}(r,e)$
    and some conjectures}}

The deformation theoretic methods used in \cite{Le Potier} and \cite{Hein} allow one
to prove results similar to Theorem \ref{10} for pluritheta linear series on $U_{X}(r,e)$
(with some extra effort due to the fact that in this case one has to control the determinant
of the complementary vector bundle). Since the method used in this paper is of a different
nature, a generalization along those lines is not immediately apparent. Instead we propose
the formalism of Verlinde bundles, which we developed in \cite{Popa2}. This comes with the
advantage that it applies automatically as soon as one has results on $SU_{X}(r,e)$ and also
suggests what the optimal bounds should be. Moreover, the method equally applies to other
linear series on $U_{X}(r,e)$, as we will see shortly.

Fix a generic vector bundle $F\in U_{X}(r_{1},r_{1}(g-1)-e_{1})$, where as usual $r_{1}=r/h$ and $e_{1}=e/h$.
To it we can associate the theta divisor $\Theta_{F}$ on $U_{X}(r,e)$
supported on the set
$$\Theta_{F}=\{E\in U_{X}(r,e)|~h^{0}(E\otimes F)\neq 0\}.$$
Fix also $L\in {\rm Pic}^{e}(X)$.
The $(r,e,k)$-\emph{Verlinde bundle} $E_{r,e,k}$ associated to these choices is defined as (cf. \cite{Popa2} \S6):
$$E_{r,e,k}(=E_{r,e,k}^{F,L}):={\pi_{L}}_{*}\OO_{U}(k\Theta_{F}),$$
where $\pi_{L}$ is the composition:
$$\pi_{L}:U_{X}(r,e)\overset{{\rm det}}{\longrightarrow} {\rm Pic}^{e}(X)
\overset{\otimes L^{-1}}{\longrightarrow} J(X).$$
This is a vector bundle on $J(X)$ of rank equal to the Verlinde number $h^{0}(SU_{X}(r,e),\LL^{k})$.
The following results are proved in \cite{Popa2}:

\begin{theorem}(\cite{Popa2} 6.4 and 5.3)\label{13}
$\OO_{U}(k\Theta_{F})$ is globally generated on $U_{X}(r,e)$ as long as $\LL^{k}$
is globally generated on $SU_{X}(r,e)$ and $E_{r,e,k}$ is globally generated.
Moreover, $\OO_{U}(k\Theta_{F})$ is not globally generated for $k\leq h$.
\end{theorem}
 
\begin{proposition}(\cite{Popa2} 5.2)\label{14}
$E_{r,e,k}$ is globally generated for $k\geq h+1$ and this bound is optimal.
\end{proposition}
We immediately obtain by using \ref{10} the following bound, where $s$ is the invariant defined in
the previous section:

\begin{theorem}\label{12}
$\OO_{U}(k\Theta_{F})$ is globally generated on $U_{X}(r,e)$ for $k\geq {\rm max}\{\frac{(r+1)^{2}}{4r} h,
\frac{r^{2}}{4s} h\}$.
\end{theorem}
In fact the theorem is a special case of the more general \ref{11} that we will treat at the end of the section.
Right now it is interesting to see how these bounds relate to possible optimal bounds and discuss
some conjectures and questions in this direction. Given the shape of the result, we will carry out this discussion in the
case of $SU_{X}(r)$ and $SU_{X}(r,\pm 1)$, based on the results \ref{30} and \ref{31}. 
A similar analysis can be applied
to any other case, but we will not give any details here.

We begin by looking at degree $0$ vector bundles, where global generation is attained for $k\geq
\frac{(r+1)^{2}}{4}$, with the improvement \ref{6} in the case of rank $2$, when $k\geq 1$ suffices.
In view of \ref{13}, the bound in \ref{12} is optimal in the case of rank $2$ and rank $3$ vector bundles.
 
\begin{corollary}\label{50}
Let $N\in {\rm Pic}^{g-1}(X)$. Then:
\newline
(i) $\OO_{U}(3\Theta_{N})$ is globally generated on $U_{X}(2,0)$.
\newline
(ii) $\OO_{U}(4\Theta_{N})$ is globally generated on $U_{X}(3,0)$.
\end{corollary}
This could be seen as a natural extension of the classical fact that $\OO_{J}(2\Theta_{N})$
is globally generated on $J(X)\cong U_{X}(1,0)$. In presence of this evidence
it is natural to conjecture that this is indeed the case for any rank:
 
\begin{conjecture}\label{15}
For any $r\geq 1$, $\OO_{U}(k\Theta_{N})$ is globally generated on $ U_{X}(r,0)$ for
$k\geq r+1$.
\end{conjecture}
This is the best that one can hope for and there is a possibility that it might be
a little too optimistic, or in other words that Corollary \ref{50} might be an accident of
low values of a quadratic function. On the other
hand if that is the case, the theorem should be very close to being optimal. Turning to $SU_{X}(r)$,
in \cite{Popa} \S3 we showed that, granting the strange duality conjecture, the optimal bound
for the global generation of $\LL^{k}$ should also go up as we increase the rank $r$. The underlying
reason (without specifying the actual numbers) is the following: assume that we are given a vector
bundle $E$ such that $h^{0}(E\otimes \xi)\neq 0$ for all $\xi\in {\rm Pic}^{0}(X)$ (for examples see
\cite{Raynaud}, \cite{Popa} or \cite{Popa2}). If we choose some complementary bundle $F$ (i.e.
$\chi(E\otimes F)=0$), of rank $t$, then a theorem of Lange and Mukai-Sakai (see \cite{Lange} and
\cite{Mukai}) asserts that $F$ admits a line subbundle of degree $\geq \mu(F)-g+g/t=g/t-1-\mu(E)$.
For small $t$ (with respect to $r$), in most examples mentioned above it happens that this number is
positive, which automatically implies that $h^{0}(E\otimes F)\neq 0$ for all such $F$. This prevents
the global generation of a certain multiple of $\LL$ depending on the rank of $F$. The case of rank $2$
vector bundles \ref{6} suggests though that we could ask for a slightly better result than for
$U_{X}(r,0)$, but unfortunately further evidence is still missing:
 
\begin{conjecture/question}
Is $\LL^{k}$ globally generated on $SU_{X}(r)$ for $k\geq r-1$?
\end{conjecture/question}
Note also that in view of \ref{13} and \ref{14} any positive answer in the range
$\{r-1,r,r+1\}$ would imply the optimal conjecture \ref{15}.
 
In the case of $SU_{X}(r,\pm 1)$ \ref{31} and \ref{13} give that $\OO_{U}(k\Theta_{F})$ is globally
generated for $k\geq {\rm max}\{2, r-1\}$, while $\OO_{U}(\Theta_{F})$ cannot be. We obtain thus again
optimal bounds for rank $2$ and rank $3$ vector bundles.
 
\begin{corollary}
$\OO_{U}(2\Theta_{F})$ is globally generated on $U_{X}(2, 1)$ and $U_{X}(3, \pm 1)$.
\end{corollary}
Note also that for all the examples of special vector bundles constructed in \cite{Raynaud},
\cite{Popa} and \cite{Popa2} we have $h\neq 1$, therefore theoretically an optimal bound that
does not depend on the rank $r$ is still possible. It is natural to ask if the best possible
result always holds:
 
\begin{question}
Is $\OO_{U}(2\Theta_{F})$ on $U_{X}(r,\pm 1)$, and so also $\LL^{2}$ on $SU_{X}(r,\pm 1)$,
globally generated? More generally, is this true whenever $r$ and $d$ are coprime?
\end{question}

We conclude the section with a generalization of Theorem \ref{12}.
For simplicity we present it only in the degree $0$ case, but the extension to other degrees 
is immediate. Recall from \cite{Drezet} Theorem C
that for $N\in {\rm Pic}^{g-1}(X)$, ${\rm Pic}(U_{X}(r,0))\cong \ZZ\cdot \OO(\Theta_{N})\oplus
\dt^{*}{\rm Pic}(J(X))$.
The method provided by the Verlinde bundles allows one to study effective global generation
for ``mixed'' line bundles of the form $\OO(k\Theta_{N})\otimes \dt^{*}L$ with
$L\in {\rm Pic}(J(X))$. Concretely we have the following cohomological criterion (assume $r\geq 2$): 

\begin{theorem}\label{11}
$\OO(k\Theta_{N})\otimes {\rm det}^{*}L$ is globally generated if $k\geq \frac{(r+1)^{2}}{4}$ and
$$h^{i}(\OO_{J}((kr-r^{2})\Theta_{N})\otimes L^{\otimes r^{2}}\otimes \alpha)=0,~\forall
i>0,~\forall \alpha \in {\rm Pic}^{0}(J(X)).$$   
\end{theorem}
\begin{proof}
By the projection formula, for every $i>0$ we have:
$$R^{i}\dt_{*}(\OO_{U}(k\Theta_{N})\otimes \dt ^{*}L)\cong R^{i}\dt _{*}\OO_{U}(k\Theta_{N})\otimes
L=0.$$
Also the restriction of $\OO_{U}(k\Theta_{N})\otimes \dt ^{*}L$ to any fiber $SU_{X}(r,\xi)$ of
the determinant map is isomorphic to $\LL^{k}$ and so globally generated for
$k\geq \frac{(r+1)^{2}}{4}$.
It is a simple consequence of general machinery, described in
\cite{Popa2} (7.1), that in these conditions the statement
holds as soon as
$$\dt_{*} (\OO_{U}(k\Theta_{N})\otimes \dt^{*} L)\cong E_{r,k}\otimes L$$
is globally generated on $J(X)$, where $E_{r,k}$ is a simplified notation for $E_{r,0,k}$.
To study this we make use, as in \cite{Popa2}, of a
cohomological criterion for global generation of vector bundles on abelian varieties
due to Pareschi \cite{Pareschi} (2.1). In our particular setting it says that
$E_{r,k}\otimes L$ is globally generated if there exists some ample line bundle $A$ on $J(X)$ such that
$$h^{i}(E_{r,k}\otimes L\otimes A^{-1}\otimes \alpha)=0, ~\forall i>0,
~\forall \alpha\in {\rm Pic}^{0}(J(X)).$$
We chose $A$ to be $\OO_{J}(\Theta_{N})$, where $\Theta_{N}$ is the theta divisor on $J(X)$
associated to $N$. The cohomology vanishing that we need is true if it
holds for the pullback of $E_{r,k}\otimes L\otimes
\OO_{J}(-\Theta_{N})\otimes \alpha$ by any finite cover of $J(X)$.
But recall from \cite{Popa2} (2.3) that $r_{J}^{*}E_{r,k}\cong \bigoplus \OO_{J}(kr\Theta_{N})$,
where $r_{J}$ is the multiplication by $r$.
Since $r_{J}^{*}L\equiv L^{\otimes r^{2}}$, via pulling back by $r_{J}$ the required
vanishing certainly holds if
$$h^{i}(\OO_{J}((kr-r^{2})\Theta_{N})\otimes L^{\otimes r^{2}}\otimes \alpha)=0,~\forall
i>0,~\forall \alpha \in {\rm Pic}^{0}(J(X)).$$
\end{proof}

\begin{corollary}
If $l\in \ZZ$, $\OO(k\Theta_{N})\otimes {\rm det}^{*}\OO_{J}(l\Theta_{N})$ is globally generated for
$$k\geq {\rm max}\{r+1-lr, \frac{(r+1)^{2}}{4}\}.$$
\end{corollary}

\section{\textbf{An application to surfaces \`a la Le Potier}}

Another, in some sense algorithmic, application of the effective bound
\ref{10} can be given following the 
paper of Le Potier \cite{Le Potier}. By a simple use of a restriction theorem due to Flenner
\cite{Flenner} (cf. also \cite{Le Potier2} \S11),
Le Potier shows that effective results for the determinant bundle $\LL$ induce effective results
for the Donaldson determinant line bundles on moduli spaces of semistable sheaves on surfaces.
For the appropriate definitions and basic results, the reader can consult \cite{Lehn} \S8. 
Using the uniform bound $k\geq (r+1)^{2}/4$ that works on every moduli space $SU_{X}(r,d)$, 
the result can be formulated as follows:

\begin{theorem}\label{60}
Let $(X,\OO_{X}(1))$ be a polarized smooth projective surface and $L$ a line bundle on $X$. 
Let $M=M_{X}(r,L,c_{2})$ be the 
moduli space of semistable sheaves of rank $r$, fixed determinant 
$L$ and second Chern class $c_{2}$ on $X$ and denote $n={\rm deg}(X)=\OO_{X}(1)^{2}$ and 
$d=n[\frac{r^{2}}{2}]$. If $\mathcal{D}$ is the Donaldson determinant line bundle
on $M$, then $\mathcal{D}^{\otimes p}$ is globally generated for $p\geq d\cdot \frac{(r+1)^{2}}{4}$
divisible by $d$.
\end{theorem} 
Note that it is not true that $\mathcal{D}$ is ample, which accounts for the formulation of the 
theorem. The significance of the map to projective space given by some multiple of $\mathcal{D}$ 
is well known. Its image is the moduli space of $\mu$-semistable sheaves
and in the rank $2$ and degree $0$ case this is homeomorphic to the 
Donaldson-Uhlenbeck compactification of the moduli space 
of $ASD$-connections in gauge theory, the map 
realizing the transition between the Gieseker and Uhlenbeck points of view (see e.g. \cite{Lehn} 
\S8.2 for a survey). Better bounds for 
the global generation of the Donaldson line bundle thus limit the dimension of an ambient projective  
space for this moduli space. The main improvement brought by the results in the 
present paper comes from the fact that our result is not influenced by the genus of the curve 
given by Flenner's theorem. Effectively that reduces the bound given in \cite{Le Potier} \S3.2 by an 
order of four, namely from a polynomial of degree $8$ in the rank $r$ to a polynomial of 
degree $4$.  
\medskip
\newline
\emph{Sketch of proof of \ref{60}.}(cf. \cite{Le Potier} 3.6)
In analogy with the curve situation, given $E\in M$, the problem is to find a complementary 
$1$-dimensional sheaf $F$ on $X$ such that $h^{1}(E\otimes F)=0$. Flenner's theorem says that 
there exists a smooth curve $C$ belonging to the linear series $|\OO_{X}(d)|$ such that $E_{|C}$ 
is semistable. By theorem \ref{10}, on $C$ one can find for any $k\geq (r+1)^{2}/4$ a vector bundle
$V$ of rank $kr_{1}$ such that $h^{1}(E\otimes V)=0$. The $F$ that we are looking for is obtained 
by considering $V$ as a $1$-dimensional sheaf on $X$ and a simple computation shows that if 
$p=dk$, this gives the global generation of $\mathcal{D}^{\otimes p}$.

\begin{remark}
Depending on the values of the invariants involved, this bound may sometimes be improved to a 
polynomial of degree $3$ in $r$, according to the precise statement of Theorem \ref{10}.
\end{remark}

\end{document}